\documentclass[10pt]{article}
\usepackage{amsmath,amsfonts,amsthm,amssymb,amscd}

\catcode`\é=\active\def é{\'e} \catcode`\è=\active\def è{\`e}
\catcode`\ç=\active\def ç{\c{c}} \catcode`\à=\active\def à{\`a}
\catcode`\ê=\active\def ê{\^e} \catcode`\ù=\active\def ù{\`u}
\catcode`\ô=\active\def ô{\^o} \catcode`\û=\active\def û{\^u}
\catcode`\î=\active\def î{\^{\i}} \catcode`\â=\active\def â{\^a}
\catcode`\ö=\active\def ö{\"o}

\newcommand{\be}{\begin{equation}}
\newcommand{\ee}{\end{equation}}
\newcommand{\bs}{\begin{split}}
\newcommand{\es}{\end{split}}
\newcommand{\ba}{\begin{align}}
\newcommand{\ea}{\end{align}}
\newcommand{\basl}[1]{\begin{align}\begin{split}\label{#1}}
\newcommand{\bas}{\begin{align}\begin{split}}

\newtheorem{theo}{Theorem}[section]
\newtheorem{prop}[theo]{Proposition}
\newtheorem{lemm}[theo]{Lemma}

\newtheorem{defi}[theo]{Definition}

\def\carre{\hbox{\vrule \vbox to 7pt{\hrule width 6pt \vfill \hrule}\vrule }}

\newcommand\N{\mathbb{N}}

\newcommand\R{\mathbb{R}}

\title{Beals characterization of
pseudodifferential operators\\ in Wiener spaces}

\author{L. Amour\textsuperscript{1}, R. Lascar\textsuperscript{2} and
J. Nourrigat\textsuperscript{1}}

\date{ \ \textsuperscript{1}Université de Reims\hskip 1cm \ \textsuperscript{2}Université Paris 7, Denis Diderot  }

\begin{document}

\maketitle

\begin{abstract}
\noindent
The aim of this article is to prove a Beals type characterization theorem for pseudodifferential operators in Wiener spaces. The definition of pseudodifferential operators in Wiener spaces and a Calder\'on-Vaillancourt type result appear in \cite{AJN}. The set of symbols considered here is the one of \cite{AJN}. The Weyl calculus in infinite dimension considered here emphasizes the role of the Wick bi-symbols.
\end{abstract}

\parindent=0pt

\tableofcontents

\parindent = 0 cm

\parskip 10pt
\baselineskip 15pt

\section{Statement of the main result.}\label{s1}

In quantum field theory, such as quantum electrodynamics which will be considered in a forthcoming article, the set of states of the quantized field may be chosen as a symmetrized Fock space ${\cal F}_s (H_{\bf C})$
over an Hilbert space $H$. Among the operators acting in such spaces, those coming from the Weyl calculus in infinite dimension and recently introduced in  \cite{AJN} (see also  in \cite{AJN2} the case of the large but finite dimension)  may have applications to modelling  the interaction of the quantized field with a fixed particle of spin $1/2$.
These applications will be developed in a next article, but we need some properties which are not in \cite{AJN} and that we  present it here.

 We note by  $H$ a real separable space and by $H_{\bf C}$ the complexified. The norm of $H$ is noted by
 $|\cdot |$ and the scalar product of two elements $a$ and  $b$ of $H$
 is by $a \cdot b$. The norm of an element of $H^2$ is denoted by  $|\cdot |$. For all $X = (x ,\xi)$ and $Y
= (y ,\eta)$
 in $H^2$, we set
\be\label{1.1}
X \cdot  \overline Y = (x+i \xi) \cdot (y-i\eta),
 \qquad \sigma (X , Y) = y\cdot \xi - x\cdot \eta.   \ee

We recall that ${\cal F}_s (H_{\bf C})$ is the completion of the direct sum of the subspaces ${\cal F}_n$ ($n\geq 0$) where ${\cal F}_0$ is one dimensional and represents the vacuum, while ${\cal F}_1 = H_{\bf C} $ and ${\cal F}_n$ ($n\geq 2$) is the $n-$ fold symmetrized tensor product representing the $n$ particles states. This space is not very convenient for the Weyl calculus since we have to write down integrals but it is isomorphic to some $L^2$ space on a suitable Banach $B$ endowed with a gaussian measure.

It is known that, for any  separable real Hilbert space $H$ there exists,

\hskip 1cm  - a Banach space $B$ containing $H$,

\hskip 1cm  - a  gaussian measure $\mu _{B , h }$ with variance $h$ on the $\sigma-$algebra of the Borel sets of $B$, for all $h>0$,

satisfying some assumptions  we formulate here in saying that
 $(i, H, B)$ is an abstract Wiener space
(where $i$ is the injection from  $H$ into $B$). See \cite{G1}\cite{G2}\cite{KU} and \cite{AJN} for precise conditions which should be fullfilled by $B$.  See also  \cite{G3} (example 2, p. 92) for a standard way of construction of a space $B$ satisfying the assumptions.

Identifying $H$ with its dual, one has,
\be\label{1.2} B' \subset H' = H \subset B.  \ee
If $H$ is finite dimensional, we have  $B= H$ and for all Borel sets $\Omega$ in $H$,
\be\label{1.3} \mu _{H , h} (  \Omega)  = (2\pi h)^{-{\rm dim} (E)/2}
 \int _{\Omega}   e^{-{|y|^2 \over 2h}} dy.
      \ee
In the general case, the symmetrized Fock space ${\cal F}_s (H_{\bf C})$ (\cite{SE},\cite{RS}) is isomorphic to the space $L^2(B, \mu_{B , h/2})$ (see \cite{J}\cite{SI}). The complexified $H_{\bf C}
 \subset {\cal F}_s (H_{\bf C}) $ is identified with a closed subset of $L^2(B, \mu_{B , h/2})$ which in field theory is the subspace corresponding to the states of the field with exactly one particle.

The Weyl calculus in infinite dimension of \cite{AJN} allows to associate to some suitable functions $F$ on the Hilbert space $H^2$, bounded and unbounded operators in ${\cal F}_s (H)$ (or in  $L^2(B, \mu_{B , h/2})$).  Let us first recall the assumptions filled by  functions  $F$.

\begin{defi}\label{d1.1}  Let $(i, H, B)$ be a Wiener space satisfying  (\ref{1.2}).
 We choose a Hilbert
basis $(e_j )_{(j\in \Gamma)}$ of $H$, each vector belonging to
$B'$, indexed by a countable set $\Gamma$. Set $u_j = (e_j , 0)$ and
$v_j = (0, e_j)$ $(j\in \Gamma)$. A multi-index is a map
 $(\alpha , \beta )$ from $\Gamma $ into
$\N \times \N$ such that $\alpha_j = \beta _j = 0$ excepted for a
finite number
 of indices. Let $M$ be a nonnegative real number, $m$  a nonnegative integer
 and  $\varepsilon = (\varepsilon_j )_{(j \in \Gamma)}$ a family of
nonnegative real numbers.  One denotes by $ S_m(M, \varepsilon)$ the
set of bounded continuous functions $ F:H^2\rightarrow {\bf
C}$ satisfying the following conditions. For every multi-index
 $(\alpha , \beta)$ such that $0 \leq \alpha_j \leq
m$ and $0 \leq \beta_j \leq m$ for all $j\in \Gamma$, the following
derivative,
\be\label{1.4}\partial_x^{\alpha}\partial_{\xi}^{\beta}  F =  \left [\prod _{j\in \Gamma }
\partial _{u_j} ^{\alpha_j} \partial _{v_j} ^{\beta_j}\right ]  F  \ee
is well defined, continuous on
 $H^2$ and satisfies, for every $(x , \xi)$ in  $H^2$,
\be\label{1.5}\left |\partial_x^{\alpha}\partial_{\xi}^{\beta}   F(x , \xi)
\right |    \leq M \prod _{j\in \Gamma } \varepsilon_j ^{\alpha_j +
\beta_j}. \ee
\end{defi}

For each summable sequence
 $(\varepsilon_j)$, the first step in \cite{AJN} is to associate to each function $F$ in $S_2(M, \varepsilon)$, a quadratic form $Q_h^{weyl} (F)$ on a dense subset ${\cal D}$ (see
 Definition \ref{d2.1} above), and not an operator on the above Hilbert spaces.

One may also associate a quadratic form  $Q_h^{weyl} (F)$  on ${\cal D}$ with  symbols $F$ which are not in the above set, in particular if they are not bounded.
To do it, it is sufficient that the two conditions below are satisfied:

 (H1) The function $F: H^2 \rightarrow {\bf C} $ has a stochastic extension
 $\widetilde F : B^2  \rightarrow {\bf C} $ in $L^1 ( B^2 , \mu _{B^2 , h/2})$ (see definition
 4.4 of \cite{AJN} which recall and adapt a previous definition of L. Gross \cite{G1}.

 (H2) The action on $|\widetilde F|$ of the following heat operator
\be\label{1.6} ( H_{h/2} |\widetilde F | ) (X) = \int _{B^2} |\widetilde F(X+Y) |
 d\mu _{B^2 , h/2} (Y)  \hskip 2cm X\in H^2\ee
 is polynomially bounded, i.e., it satisfies for $m\geq 0$ and
 $C>0$,
\be\label{1.7} ( H_{h/2} |\widetilde F | ) (X) \leq C (1+ |X|)^m \ee
 (that is to say that the norm in formula  (12) in \cite{AJN} is finite).

In Theorem \ref{t2.2}, we  recall the construction of $Q_h^{weyl} (F)$ in a slightly simplified way, but the construction in \cite{AJN} uses the analog in infinite dimension of Wigner functions which may have its own interest. The hypotheses 1 and 2 are satisfied if $F$ belongs to $S_2(M, \varepsilon)$, the sequence $(\varepsilon _j)$ being
  summable. Inequality (\ref{1.7}) is then satisfied with $C=M$ and $m=0$.  See others examples in Section \ref{s2}.

Next, as  shown in  \cite{AJN} (Theorem 1.4), if $F$ belongs to $S_2(M, \varepsilon)$ then $Q_h^{weyl} (F)$
is the  quadratic form of a bounded operator in  $L^2(B, \mu_{B , h/2})$ or equivalently, bounded in  ${\cal F}_s(H_{\bf C})$. In addition, this operator satisfies, if $0 < h < 1$,
\be\label{1.12}  \Vert  Op_h^{weyl} (F) \Vert \leq M  \prod _{j\in \Gamma} (1 + 81\pi h S_{\varepsilon}  \varepsilon_j ^2) \ee
where
\be\label{1.13}  S_{\varepsilon}  = \sup_{j\in \Gamma}  \max(1, \varepsilon_j ^2 ).\ee
The hypothesis (H2) in Theorem
1.4 in \cite{AJN}, which not mentioned here, is always satisfied if
$F$ belongs to $S_2(M, \varepsilon)$ and if the sequence $(\varepsilon_j )$ is summable (Proposition
8.4 in \cite{AJN}).

We have now to define and to compute, commutators of these operators with momentum and position operators. In finite dimension $n$, theirs compositions and commutators are a classically defined as operators from ${\cal S} (\R^n)$ into
${\cal S}' (\R^n)$. In our case, ${\cal S} (\R^n)$ is replaced by  space
${\cal D}$ of Definition \ref{d2.1}. In the absence of an analog of ${\cal S}' (\R^n)$,
we prefer instead to use quadratic forms on  ${\cal D}$ (see \cite{RS}). We then consider mappings  $(f , g) \rightarrow A(f , g)$ on ${\cal D} \times {\cal D}$ that are linear in $f$ and antilinear in $g$. A notion of continuity is given in Section \ref{s2}.

One may define two compositions (left and right) of a quadratic form $Q$ on the space  ${\cal D}$ of Definition \ref{d2.1}
with an operator $A : {\cal D} \rightarrow {\cal D}$ whose formal adjoint $A^{\star}$ also maps ${\cal D}$ into  ${\cal D}$.
One set, for all $f$ and $g$ in ${\cal D}$,
\be\label{a2}( Q \circ A ) (f , g) = Q ( Af , g),\qquad( A \circ Q ) (f , g) = Q ( f , A^{\star} g).\ee
One then define the commutator $[A, Q ]$ and $({\rm ad} A) Q$  as the following quadratic form,
\be\label{a3}[A, Q ] (f , g) = Q ( f , A^{\star} g )  - Q ( Af , g). \ee
Thus, one can define the iterated bracket
 $({\rm ad} A_1)  \dots ({\rm ad} A_n)  Q$ if  $A_1$,
\dots $A_n$ are operators from ${\cal D}$ into ${\cal D}$.

We  see in Proposition \ref{p2.3} that one may associate with each continuous linear form $G$ on $H^2$, not only a quadratic form $Q_h^{weyl} (G)$, but also an operator $Op_h^{weyl} (G)$ from ${\cal D}$ to ${\cal D}$. This Weyl operator is the Segal field, up to a numerical factor, and may be directly defined in ${\cal F}_s (H)$ using creation and annihilation operators, without using the Weyl calculus. In particular, when
$F(x , \xi) = a \cdot x$ with $a$ in $H$, the corresponding Weyl operator will be denoted $Q_h (a)$ (position operator). When $F(x , \xi) = b \cdot \xi $,
when $b$ in $H$, the operator will be denoted $P_h(b)$ (momentum operator).

If $F$ belongs to $S_m (M, \varepsilon)$ and $G$ is a continuous linear form on $H^2$  then  Proposition \ref{p2.6} allows us to extend the following result which is well-known in finite dimension,
\be\label{a1} [Q_h^{weyl } (F)  , Op_h^{weyl } (G)] = {h\over i} Q_h^{weyl } ( \{ F , G \} ).\ee
In particular, if $(e_j)$ is the Hilbertian basis of $H$ chosen to define our  sets of symbols then  equality (\ref{a1}) gives,
$$ [Q_h (e_j) , Q_h^{weyl } (F)  ] = - {h\over i} Q_h^{weyl } \left ( {\partial  F
\over \partial \xi_j} \right ), $$
$$ [P_h (e_j) , Q_h^{weyl } (F)  ] =  {h\over i} Q_h^{weyl } \left ( {\partial  F
\over \partial x_j} \right ). $$
One may iterate and consider iterated commutators while restricting ourselves to some set of multi-indices. We denote by ${\cal M}_m$ the set of pairs $(\alpha , \beta)$ where   $\alpha = (\alpha_j )_{(j\in \Gamma)}$ and
$\beta = (\beta_j )_{(j\in \Gamma)}$ are sequences of nonnegative integers such that $\alpha_j = \beta _j =0$ except for a finite number of indices $j$, and such that $\alpha _j \leq m$ and $\beta_j \leq m$ for all $j\in \Gamma$.
One associates to each multi-index $(\alpha , \beta)$ the following iterated commutator,
$$ ({\rm ad}P_h )^{\alpha} ({\rm ad}Q_h )^{\beta} Q_h ^{weyl} (F) = \prod _{j\in \Gamma} ( ad P_h(e_j) ) ^{\alpha _j}
\prod _{k\in \Gamma} ( ad Q_h(e_k) ) ^{\beta _k}Q_h ^{weyl} (F). $$

In the same way, if $F$ is in $S_m(M, \varepsilon)$ and if $(\alpha , \beta)$
is in ${\cal M}_p$,  $p \leq m-2$,
$$ ({\rm ad}P_h )^{\alpha} ({\rm ad}Q_h )^{\beta} Q_h ^{weyl} (F) = (-1)^{|\beta |}
(h/i) ^{|\alpha + \beta |} Q_h ^{weyl} (\partial_x^{\alpha } \partial_{\xi}^{\beta } F).  $$
From Theorem 1.4
in \cite{AJN}, the above Weyl quadratic form  is associated to a bounded operator in  $L^2(B , \mu_{B , h/2})$,  denoted as below and verifiying,
\be\label{1.15}\Vert ({\rm ad}P )^{\alpha} ({\rm ad}Q )^{\beta} Op_h^{weyl}(F)  \Vert \leq
 M  \prod _{j\in \Gamma} (1 + 81\pi h S_{\varepsilon}  \varepsilon_j ^2)
  \prod _{j\in \Gamma} (h\varepsilon _j )^{\alpha _j +\beta _j}.\ee
The purpose of this work is to prove the reciprocal statement, as Beals \cite{Bea} did in finite dimension
(see also \cite{BO1}\cite{BO2} and \cite{BO-C} for adaptations to other classes of symbols in finite dimension).

\begin{theo}\label{t1.2}  Let $(i, H, B)$ be a Wiener space  satisfying (\ref{1.1}).
Let $A_h$ be a bounded operator in $L^2(B, \mu_{B, h/2})$. Let
$(e_j)$ $(j\in \Gamma)$  a Hilbertian basis of $H$ consisting of elements in $B'$.
Let $M > 0$ and let $(\varepsilon _j ) _{(j \in \Gamma)}$  a
 summable sequence of real numbers. Let $m \geq 2$. Suppose that,
for all $(\alpha , \beta )$  in ${\cal M}_{m+4}$, the commutator
$({\rm ad}P )^{\alpha} ({\rm ad}Q )^{\beta}  A_h$ (being a
priori defined as a quadratic form on ${\cal D}$) is bounded
in $L^2(B, \mu_{B, h/2})$ and that,
\be\label{1.15b}   \Vert ({\rm ad}P )^{\alpha} ({\rm ad}Q )^{\beta}  A_h \Vert \leq
 M \prod _{j\in \Gamma} (h\varepsilon _j )^{\alpha _j +\beta _j}.\ee
Then, if $0 < h < 1$,  there exists a function $F_h $ in $S_m(M', \varepsilon)$ with,
\be\label{1.16}  M' = M \prod _{j\in \Gamma} (1 + K S_{\varepsilon} ^ 2 h \varepsilon _j ^2)  \ee
where $K$ is a universal constant, and $S_{\varepsilon}$  is defined in (\ref{1.13}),
 such that the Weyl operator $Op_h^{weyl} (F )$ associated to $F$ is equal to $A_h$.
 \end{theo}

Section \ref{s2} introduces various results concerning the  Weyl calculus in infinite dimension intended to be used in an upcoming work. Sections \ref{s3} to \ref{s7}
are devoted to  proof of Theorem \ref{t1.2}. Section \ref{s8} applies this
theorem  to  composition of two operators defined
 by the Weyl calculus. We show that the composition is also defined by this  calculus, but
we do not give any results on the possible asymptotic expansion of its symbol, this result being used in a forthcoming article.

\section{Weyl calculus in infinite dimension.}\label{s2}

\subsection{Coherent states.}\label{s2.A}

For  $X=(a , b)$ in $H^2$, and  all $h>0$, one defines $\Psi_{X ,h}$ the corresponding coherent state (\cite{Ber}\cite{C-R}\cite{F}), they belong to ${\cal F}_s (H_{\bf C})$ and are defined by,
\be\label{2.1} \Psi_{(a , b) , h}  = \sum _{n\geq 0}
 {e^{-{|a|^2+ |b|^2  \over 4h}} \over (2h)^{n/2} \sqrt {n!} } (a+ib) \otimes \cdots \otimes (a+ib).\ee

In view of the isomorphism from ${\cal F}_s (H_{\bf C})$ in $ L^2(B , \mu_{B , h/2})$,
each element $a$ of $H \subset {\cal F}_s (H_{\bf C})$ is seen as a
function in $L^2(B , \mu_{B , h/2})$ denoted $x \rightarrow \sqrt {h} \ell_a(x)$.
When $a$ is in $B' \subset H$, one has $\ell_a(x) = a(x)$. When $a$ is in $H$, it is approximated by a sequence $(a_j)$ in $B'$, we then show that the sequence $\ell_{a_j}$
is a Cauchy sequence in $ L^2(B , \mu_{B , h/2})$ and we denote  by $\ell _a$ its limit.
With the same isomorphism, the coherent state $ \Psi_{(a , b) , h}$ defined in (\ref{2.1}) becomes,
\be\label{2.2}\Psi_{X , h} (u) = e^{{1\over h} \ell _{ (a+ib)} (u)  -{1\over 2h}|a|^2 -
{i\over 2h} a\cdot b},\quad X = (a , b) \in H^2,\quad {\rm a.e.}\ u\in B. \ee
We see, for all $X = (x , \xi)$ and $Y = (y, \eta)$, with the notation (\ref{1.1}), that
\be\label{2.3}< \Psi_{X h} , \Psi _{Yh}> =e^{-{1\over 4h}(|X|^2 +|Y|^2) + {1\over 2h} X \cdot \overline Y }.\ee
In particular,
\be\label{2.4}|< \Psi_{X h} , \Psi _{Yh}>|  =e^{-{1\over 4h}|X-Y|^2 }.\ee

We call  Segal Bargmann transform (\cite{HA}) of $f$ the function
\be\label{2.5}(T_hf) (X) = { < f , \Psi_{Xh} > \over  < \Psi_{0h} , \Psi_{Xh} >},\qquad X\in H^2.\ee
We know that $T_hf$ admits a stochastic extension $\widetilde T_hf $ in $L^2(B^2 , \mu _{B^2, h})$ and we know
that,  $\widetilde T_h$ is a partial isometry from  $L^2(B , \mu _{B, h/2})$ into $L^2(B^2 , \mu _{B^2, h})$.

\subsection{The space ${\cal D}$ and  Wick symbols.}\label{s2B}

\begin{defi}\label{d2.1} For all subspaces $E$ of finite dimension in $H$,
 ${\cal D}_E$ denotes the space of functions $f : B \rightarrow {\bf C}$ such that,

i) the  function $f$ is written under the form $ \widehat f \circ P_E$, where $\widehat f $
is a continuous function from $E$ in ${\bf C}$ and $P_E$ is the mapping from  $B$ in
$E$ defined as follows, choosing an orthonormal basis
$\{ u_1 , ... u_n \}$ of $E$,
\be\label{2.6}P_E(x) = \sum _{j=1}^n \ell _{u_j} (x) u_j,\quad a.e. \  x\in B\ee
(the map $P_E$ is independent of the chosen basis).

\smallskip ii) the function $E^2 \ni X \rightarrow < f, \Psi_{X h}>$ (scalar product in
 $L^2(B , \mu_{B , h/2})$) is in the
Schwartz space ${\cal S} (E^2)$.

\smallskip We shall denote by  ${\cal D}$ the union of all spaces ${\cal D}_E$.
\end{defi}

We observe that the coherent states belong to ${\cal D}$.  The condition ii)
is equivalent to say that the function $ \widehat f $ of i) is such that the function
\be\label{2.7} E \ni u \rightarrow \widehat f (u) e^{-{|u|^2 \over 2h}}\ee
belongs to ${\cal S} (E)$. One says that a quadratic form
$Q$ on ${\cal D}$ is {\it continuous} if, for all $E \subset H$ of finite dimension, there exists $C>0$ and $m\geq 0$ such that, for all $f$ and $g$ in ${\cal D}_E$,
\be\label{2.8} |Q(f, g)| \leq C I(E, m) (f) I(E, m) (g)\ee
where
\be\label{2.9} I(E, m) (f) = \int _{E^2} |< f, \Psi_{X h}>| (1+|X|)^m dX.\ee
One says that a linear mapping $T $ in ${\cal D}$ is continuous if,
for all $E\subset H$ of finite dimension, there exists $F\subset H$ of finite dimension
such that $f\in {\cal D}_E$ implies $Tf \in {\cal D}_F$ and if, for all
integer $m$, there exists $C$ and $m'$ such that,
\be\label{2.10} I(F, m) (Tf) \leq C I(E , m') (f).\ee

We shall recall the definition of the Wick symbol and bi-symbol. If $Q$ is a quadratic form
 on  ${\cal D}$, we denote by $S_h (Q)$ the function defined on $H^2$ by,
\be\label{2.11} S_h (Q) (X , Y)=  { Q( \Psi_{X , h} ,    \Psi_{Y , h})  \over
<  \Psi_{X , h} , \Psi_{Y , h} >}.\ee
If $Q(f , g) = < Af, g>$, where $A$ is an bounded operator  in
the  Fock space ${\cal F}_s (H_{\bf C})$, or equivalently in $L^2(B, \mu_{B , h/2})$, then
the symbol $S_h (Q)$ will be also denoted $S_h (A)$.
 Let us recall that, if  $X = (x , \xi)$ is identified with $x+i \xi$, then the function $S_h(A)$ is Gateaux
holomorphic in $X$ and antiholomorpic in  $Y$.

We denote by  $\sigma_h^{wick} (Q)$ the restriction to the  diagonal of the above function,
\be\label{2.12}\sigma_h^{wick} (Q) (X) = Q( \Psi_{X , h} ,    \Psi_{X , h}).\ee

\subsection{Definition of the Weyl calculus in infinite dimension.}\label{s2.C}

If  $H= B =\R^n$ and if, say, $F$ is  a $C^{\infty }$ function
 on $\R^{2n}$ bounded together with all its derivatives, one associates with $F$ an operator $Op_h^{weyl} (F)$  satisfying,
\be\label{2.13} S_h ( Op_h^{weyl} (F) )  (X, Y) = \int _{\R^{2n}}F(Z)
e^{{1 \over  h}( X \cdot \overline Z + \overline Y \cdot Z - X \cdot \overline Y)  } d\mu _{h/2} (Z)
\ee
$$=e^{{1\over 4h}|X - Y|^2}  \int _{\R^{2n}} F\left (Z + {X+Y\over 2}
\right ) e^{{i\over 2h}( (\xi - \eta ) \cdot  z - (x- y) \cdot \zeta  )}  d\mu _{\R^{2n}  , h/2} (Z).
$$
This equality  is proved in Unterberger \cite{U} and we use it for an extension to the infinite dimensional spaces.

The first issue is that,  the function $F$ is defined on $H^2$ according the Definition \ref{d1.1}, and
giving a meaning in infinite dimension to an integral such as the one in  (\ref{2.13}),
we have to integrate over $B^2$, where $(i, H, B)$ is a  Wiener space. Indeed,
 in infinite dimension, $H^2$ cannot be endowed with a gaussian measure which corresponds to its own norm.

We  have to be able to extend the function $F$, defined on $H^2$, to a function $\widetilde F$
defined on $B^2$. In general it is not a density extension but a type of
extension introduced  by L. Gross and named {\it  stochastic extension}. It may be found in
\cite{AJN} (Definition 4.4) where we recall a definition of this notion adapted to our purposes.
From Proposition 8.4 of \cite{AJN}, we know  that each function $F$ in  $S_1 (M, \varepsilon)$ admits
 a stochastic extension $\widetilde F$ in
$L^1 (B^2 , \mu _{B^2 , h/2})$ at least if the sequence $(\varepsilon_j)$ is summable. Moreover, the proof of Proposition 8.4 of \cite{AJN} shows that
any linear form $F$ on $H^2$   has a stochastic extension $\widetilde F$ in
$L^1 (B^2 , \mu _{B^2 , h/2})$.

By analogy with (\ref{2.13}), one expect to associate with each function $F$
satisfying the hypotheses (H1) and (H2) of Section \ref{s1},  a quadratic form $Q_h^{weyl} (F)$ on ${\cal D}$,
with  bi-symbol $S_h ( Q_h^{weyl} (F)) $ of form,
\be\label{2.14}\Phi  (X , Y)=e^{{1\over 4h}|X - Y|^2}  \int _{B^2}\widetilde F\left (Z + {X+Y\over 2}
\right ) e^{{i\over 2h}( \ell _{\xi - \eta} (z) -\ell _{ x- y} (\zeta ) )}  d\mu _{B^2 , h/2} (Z).\ee

\begin{theo}\label{t2.2} Let $F: H^2 \rightarrow {\bf C}$ be a function
satisfying the hypotheses (H1) and (H2) of  Section \ref{s1} with $m\geq 0$.
 Let $\widetilde F$ be the stochastic extension
of $F$ in $L^1 (B^2 , \mu _{B^2 , h/2})$. Then,

i)  The integral (\ref{2.14}) converges and verifies,
\be\label{2.15} |\Phi  (X , Y)|  \leq C   e^{{1\over 4h}|X-Y|^2} \left ( 1+ { |X+Y| \over 2} \right )^m.\ee
 In addition, this function is   Gateaux  holomorphic in $X$ and anti-holomorphic in $Y$.

ii) There is a continuous quadratic form $Q_h^{weyl} (F)$  on ${\cal D}$
such that $S_h ( Q_h^{weyl} (F)) = \Phi$, i.e.,
\be\label{2.16} S_h ( Q_h^{weyl} (F)) (X , Y) = e^{{1\over 4h}|X - Y|^2}  \int _{B^2}\widetilde F\left (Z + {X+Y\over 2}
\right ) e^{{i\over 2h}( \ell _{\xi - \eta} (z) -\ell _{ x- y} (\zeta ) )}  d\mu _{B^2 , h/2} (Z).\ee
\end{theo}

{\it Proof.} i) The  convergence of the integral (\ref{2.14}) and the estimate (\ref{2.15})
follow from hypothesis (H2). By a change of
variables (cf \cite{AJN}\cite{KU}), the function $\Phi$ may be also written as,
\be\label{2.17}\Phi  (X , Y)= \int _{B^2} \widetilde  F(Z)
e^{{1 \over  h}( \ell _X ( \overline Z)  + \ell _{\overline Y}  ( Z ) - X \cdot \overline Y)  } d\mu _{h/2} (Z). \ee
We deduce that it is holomorphic in $X$ and anti-holomorphic in $Y$.

ii) For all $f$ and $g$ in ${\cal D}_E$, where  $E\subset H$
is a subspace of finite dimension, set
\be\label{2.18} Q (f , g)  = \int_{E^4} \Phi (X , Y) e^{ {1\over 2h} X\cdot \overline Y }
(T_hf) (X)  \overline {(T_hg) (Y)} d\mu _{E^4 , h}(X , Y).\ee
Using (\ref{2.15}) we see that, for all $f$ and $g$ in ${\cal D}_E$,
$$ | Q(f , g) | \leq C (2\pi h)^{-2 {\rm dim} E} \int_{E^4} |< f , \Psi_{X , h}>| |< g , \Psi_{Y , h}>|
(1+|X|) (1+|Y|) d\lambda  (X , Y)$$
where $\lambda $ is the  Lebesgue measure. Consequently, for all $f$ in ${\cal D}_E$, the integral defining
$Q(f, g)$ converges. When $f$ and $g$ belongs to ${\cal D}_E$, they also are
in ${\cal D}_F$, for all subspace $F$ containing $E$. If $F$ contains $E$, then
we denote by  $S$ the orthogonal set to $F$ in $E$, and $(X_E, X_S)$ the variable of $F^2$. The transform
$T_h f$ is a function on $F^2$, independent of the variable $X_S$. We remark  that,
$$\int_{S^4} \Phi (X_E + X_S , Y_E + Y_S)e^{ {1\over 2h} X_S\cdot \overline Y_S }
d\mu _{S^4 , h}(X_S , Y_S) = \Phi (X_E  , Y_E ).$$
Indeed, the function in the  integral is holomorphic in $X_S$, anti-holomorphic in $Y_S$, and its
integral is equal to its value at $X_S = Y_S = 0$. Consequently the definition of $Q(f, g)$ is indeed coherent, whether  that $f$ and $g$
are seen as functions  in ${\cal D}_E$  or in  ${\cal D}_F$.
Let us show that the bi-symbol of $Q$ is $\Phi$. We have, for all $X= (x , \xi)$ and
$Y = (y, \eta)$ in $H^2$, if $E$ is the subspace spanned  by $x$, $\xi$, $y$
and $\eta$,
$$  {Q (\Psi_{X h} , \Psi _{Yh} ) \over < \Psi_{X h} , \Psi _{Yh}> } = \int_{E^4}
 \Phi (U , V)
{\cal B}_h(X , Y, U, V) d\mu _{E^4 , h}(U , V)$$
where ${\cal B}_h$ is a kind of reproducing kernel,
\be\label{2.19}{\cal B}_h(X , Y, U, V) = e^{ {1\over 2h} (X \cdot \overline U +  U\cdot \overline V
+   V \cdot \overline Y
 -  X\cdot \overline Y  )}.\ee
In a standard way, we have, if $\Phi$ is holomorphic in $X$, anti-holomorphic in $Y$,
\be\label{2.20} \int_{E^4}  \Phi (U , V)
{\cal B}_h(X , Y, U, V) d\mu _{E^4 , h}(U , V) =  \Phi (X , Y).\ee
It suffice to make the change of variables $U = X +S$, $V = Y + T$, and to apply the mean formula.
We then deduce  that the bi-symbol of $Q$ is indeed $\Phi$.

\hfill \carre

When $F$ belongs to $S_2 (M, \varepsilon)$,
where the sequence $(\varepsilon_j)$ is summable, we have proved in \cite{AJN} that the quadratic form $Q_h^{weyl} (F)$ is associated with a bounded operator.

\subsection{Weyl symbol and Wick symbol.}\label{s2.D}

It is sufficient to restrict equality (\ref{2.16}) to the diagonal $Y= X$ to see that,
\be\label{2.21}\sigma_h ^{wick} ( Q_h^{weyl} (F))(X) =  \int _{B^2}\widetilde F(Z + X )
 d\mu _{B^2 , h/2} (Z).\ee
For all $t>0$, the operator
\be\label{2.22}(H_tF)  (X) =
 \int _{B^2}  \widetilde F( X +Y ) d\mu _{B^2, t} (Y)\ee
is considered as the heat operator.
In the above and below integrals on $B^2$, $\widetilde F(X+Y)$ denotes the stochastic extension on $B^2$ of $H^2 \ni Y \rightarrow F(X+Y)$ for each $X$ in $H^2$, which exists since it satisfies the same hypotheses as $F$.
We then can write,
\be\label{2.23} \sigma_h ^{wick} ( Q_h^{weyl} (F)) =H_{h/2}F.\ee
Equality (\ref{2.23}) extends the standard fact in finite dimension, that the Wick symbol is obtained from the
 Weyl symbol by the action of the heat operator.
From Kuo \cite{KU} (Theorem 6.2) or Gross \cite{G4} (Proposition 9), the function $H_tF$ is continuous on $H^2$. If $H$ is of finite dimension,
we have $B= H$, $\widetilde F =F$, and $H_t F = e^{ (t/2) \Delta}  F$. Note that,
\be\label{2.24} \sup _{X \in H^2} |(H_tF)  (X)| \leq  \sup _{Z \in B^2}| \widetilde F( Z )| =
\sup _{X \in H^2} |F(X)|.\ee

\begin{prop}\label{p2.3}  If $F$ is in $S_4(M, \varepsilon)$ with some chosen basis $(e_j)$ and if the sequence $(\varepsilon_j)$ is summable, then there exists $C>0$
such that,  for all $X$ in $H^2$ and $t$ in $(0, 1)$,
\be\label{2.25} |(H_tF) (X) -  F(X)  | \leq C t.\ee
\end{prop}

 {\it Proof.}
Let $E_m$ be the subspace spanned by the $e_j$ ($j\leq m$).  We apply (\ref{2.24}) to the function
$F_m = F - F \circ \pi _{E_m}$.
We obtain, for all $X$ in $H^2$,
$$ \int _{B^2 } |( F  \circ  P _{E_m} )  (X+Y )) -  (\widetilde F_t (X+ Y )) | d\mu _{B^2, t}  (Y)
\leq   \Vert  F - F \circ \pi _{E_m} \Vert _{\infty}$$
where $\pi _{E_m} : H^2 \rightarrow E_m^2$ is the orthogonal projection and $P _{E_m} : B^2 \rightarrow E_m^2$
is its stochastic extension, defined as in (\ref{2.6}).
If $F$ is in $S_1 (M, \varepsilon)$, we have,
\be\label{2.26}\Vert  F - F \circ \pi _{E_m} \Vert _{\infty} \leq 2M  \sum _{j=p}^{\infty} \varepsilon_j.\ee
For all $m>0$ and
for all $X$ in $H^2$, we have,
$$  \int _{B^2}  F (  P _{E_m}  (X+Y ))  d\mu _{B^2, t} (Y) =
 \int _{E_m^2}  F (  (\pi _{E_m}  X)  +Y )  d\mu _{E_m^2, t} (Y).   $$
According to standard results in finite dimension, we have for all $a$ in $E_m^2$,
$$ \left |   \int _{E_m^2}  F ( a  +Y )  d\mu _{E_m^2, t} (Y) - F(a)
\right | \leq t \Vert \Delta _m  F \Vert _{\infty}$$
where
$$ \Delta_m = \sum _{j=1}^m \left ( {\partial ^2 \over \partial _{x_j}^2} +
{\partial ^2 \over \partial _{\xi_j}^2 } \right ).$$
We apply  this inequality to  $a =  \pi _{E_m}  (X)$ using again (\ref{2.26}).
Consequently, for all $t\in (0, 1)$ and $m\geq 1$,
$$ |(H_tF) (X) -  F(X) | \leq 2M  t\sum _{j=1}^m \varepsilon_j^2  + 4M  \sum _{m+1}^{\infty } \varepsilon_j.$$
 We deduce (\ref{2.25}) when $m$ goes to infinity.

 \hfill\carre

\subsection{Operators with linear symbol. Composition.}\label{s2.E}

\begin{prop}\label{p2.4}  Let $F$ be a continuous linear form  on $H^2$.  Let $Q_h^{weyl}(F)$ be the quadratic form on ${\cal D}$
defined in Theorem \ref{t2.2}. Then, there exists  an operator denoted
$Op_h ^{weyl} (F)$ from ${\cal D}$ into itself, such that
\be\label{2.27}   Q_h^{weyl} ( F) (f , g) = < Op_h ^{weyl} (F) f, g>,\qquad
(f , g)\in {\cal D}^2.\ee
\end{prop}

{\it Proof.} Let $f$ be in ${\cal D}_E$,
where $E \subset H$ is of finite dimension.  As in  Definition \ref{d2.1}, we may write, $f = \widehat f \circ P_E$,
where the function (\ref{2.7}) is in ${\cal S} (E)$.
Let  $a$ and $b$ in $H$ be such that
 $F(x , \xi ) = a \cdot x + b \cdot \xi$. Let $E_1$ be the subspace spanned
 by $E$, $a$ and $b$. Set $f_1 : E_1 \rightarrow {\bf C}$ the function
 defined by,
 $$ f_1(u) = (a+i b) \cdot u \widehat f(\pi (u)) + {h\over i}
 (\pi (b) \cdot \nabla \widehat f) (\pi (u)),\qquad u\in E_1$$
where $\pi : E_1 \rightarrow E$ is the  orthogonal projection.
We have  $OP_h^{weyl} (F) f = f_1 \circ P_{E_1}$ and this function is in
${\cal D}_{E_1}$. Thus, if $F$ is linear, the quadratic form
$Q_h^{weyl} (F)$ is associated with a continuous operator $Op_h^{weyl}(F)$
from ${\cal D}$ into ${\cal D}$.  The set of linear functions is invariant by the operator $H_{h/2}$. Consequently, the Wick symbol of $Q_h^{weyl} (F)$ is also $F$. We may write
$F(x , \xi) = P(X ) + Q(\overline X)$. Then, the bi-symbol of
$Q_h^{weyl} (F)$ is $P(X) + Q(\overline Y)$. We have, for all  $f$
 in ${\cal D}_E$, for all $Y \in (E_1)^2$,
 $$ < Op_h^{weyl} (F) f, \Psi_{Yh} > = (2\pi h) ^{-n}
 \int _{E^2} <f, \Psi_{Xh} > < Op_h^{weyl} (F) \Psi_{Xh} , \Psi_{Yh} >
dX  $$
$$ = (2\pi h) ^{-n}
 \int _{E^2} <f, \Psi_{Xh} > [P(X) + Q(\overline Y)]
 <\Psi_{Xh} , \Psi_{Yh} > dX.$$
Consequently, for all integer $m$,
$$ (1+|Y|)^m |  < Op_h^{weyl} (F) f, \Psi_{Yh} > | \hskip 8cm$$
$$ \hskip 2cm \leq C(E , E_1, h)
 \int _{E^2} (1+|X|)^{m+1}   | <f, \Psi_{Xh} > | (1+|X-Y|)^{m+1}
e^{-{1\over 4h} |Y-X|^2} dX. $$
Therefore,
$$ I(E_1 , m) (Op_h^{weyl} (F) f) \leq C(E , E_1, m, h) I(E , m+1) (f)$$
which proves the continuity of $Op_h^{weyl} (F) $ in ${\cal D}$.

\hfill\carre

Let $A$ be a continuous quadratic form  on ${\cal D}$.
Let $B : {\cal D} \rightarrow {\cal D}$ be a continuous linear mapping
with a linear Wick symbol. We recall that the quadratic forms
$  A \circ B $, $ B \circ A $ and $[A , B]$ are defined in (\ref{a2})
and (\ref{a3}).

\begin{theo}\label{t2.5}   Let $A_h$ be a bounded operator in $L^2(B , \mu _{B , h/2})$,
and set $L_h $ an operator from ${\cal D}$ into ${\cal D}$ with a Wick symbol
being a linear form $L(x , \xi) $ on $H^2$.
Let $A_h \circ B_h $ be the quadratic form  on ${\cal D}$ of their composition defined as
in Section \ref{s1}.  Then, we have,
$$ \sigma_h ^{wick} (A_h \circ L_h) = \sigma_h ^{wick} (A_h) \sigma_h ^{wick} (L_h)+ \hskip 4cm$$
$$ \hskip 4cm{h\over 2} \sum _{j\in \Gamma}  \left ( {\partial \over \partial x_j} - i
{\partial \over \partial \xi_j}\right )
\sigma_h^{wick} (A_h)
 \left ( {\partial \over \partial x_j} + i
{\partial \over \partial \xi_j}\right )
\sigma_h^{wick} (L_h) $$
This result is valid when exchanging the roles of $A_h$ and $L_h$.
\end{theo}

{\it Proof.} Set $L(x , \xi ) = a \cdot x + b \cdot \xi$ with $a$ and $b$ in $H$. Let $X $ be in
$H^2$. There exists an unitary operator $W_{X , h}$ such that $\Psi _{X , h} = W _{X , h}  \Psi _{0 , h} $.
We have,
$$ \sigma_h ^{wick} (A_h \circ L_h) (X) = < L_h \Psi_{X , h}, A_h^{\star} \Psi_{X , h } > = < f, g> $$
with $f= W_{X , h} ^{\star} L_h W_{X , h} \Psi_{0 , h} $ and $g= W_{X , h} ^{\star} A_h^{\star} W_{X , h} \Psi_{0 , h} $.
Let $T_hf$ and $T_hg $ be the Segal Bargmann transforms of $f$ and $g$ defined in (\ref{2.5}), $\widetilde T_hf $ and
$\widetilde T_hg$ being their stochastic extensions in $L^2(B^2 , \mu _{B^2, h})$. Since $\widetilde T_h$ is a partial isometry from  $L^2(B , \mu _{B, h/2})$ into $L^2(B^2 , \mu _{B^2, h})$, we have
$$ \sigma_h ^{wick} (A_h \circ L_h) (X) = \int _{B^2} \widetilde T_h f(Z) \overline { \widetilde T_h g (Z)}
d \mu _{B^2, h} (Z).$$
We also have,
$$\widetilde T_h f(Z) = L(X) + \ell _{a+ib} (z-i \zeta).$$
Since $T_hg$ is antiholomorphic then the mean formula gives,
 $$ \int _{B^2}\overline { \widetilde T_h g (Z)} d \mu _{B^2, h} (Z)= \overline {T_hg} (0) = < \Psi_{0, h}, g> =
 \sigma _h^{wick} (A_h ) (X). $$
Similarly, integrating by parts (see  Theorem 6.2 of Kuo \cite{KU}),
 for all $\gamma $ in the complexified of $H$,
 $$ \int _{B^2}\ell _{\gamma } ( z - i \zeta) \overline { \widetilde T_h g (Z)} d \mu _{B^2, h} (Z)=
h \gamma \cdot ( \partial_z - i \partial_{\zeta} ) \overline {T_hg} (0) =
 h \gamma \cdot ( \partial_x - i \partial_{\xi } )  \sigma _h^{wick} (A_h ) (X). $$
The proof of Theorem then follows.

\hfill\carre

\begin{prop}\label{p2.6}   Let $F$ be a function  in
$S_2(M, \varepsilon)$  where the sequence $(\varepsilon _j) $ is summable and let $L$ be a continuous linear form on $H^2$. Let
\be\label{b1} \Phi  = FL + {h\over 2i} \{ F , L \},\qquad\Psi  = FL - {h\over 2i} \{ F , L \}.  \ee
Then,

i) The  functions $\Phi $ and $\Psi$ satisfy hypotheses (H1) and (H2) in Section \ref{s1}

\smallskip

ii) The corresponding
Weyl forms using the Theorem \ref{t2.2} satisfy, for all $f$ and $g$ in ${\cal D}$,
$$  Q_h^{weyl } (\Phi ) (f , g) =  <   Op_h^{weyl } (L) f ,Op_h^{weyl } (F)  ^{\star} g >, $$
$$  Q_h^{weyl } (\Psi ) (f , g) =  <   Op_h^{weyl } (F) f ,Op_h^{weyl } (L)  ^{\star} g >.$$
\end{prop}

{\it Proof. i) }   Using the linearity of $G$ and the estimates
$\int_B |\ell_a(X)| |\ell_b(X)| d\mu_{B,h/2}(X) \leq C |a| |b|$, the existence of
$L^1$ stochastic extensions are obtained similarly as in the proof of the Proposition 8.4 of \cite{AJN}.
The polynomial estimate on the semigroup uses that the stochastic extension
of $X \rightarrow F(X) a.X$ is $\widetilde{F} \ell_a$ with $\int_B |\ell_a(X)| d\mu_{B,h/2}(X) \leq C |a|$.

 {\it ii) } We may write $L(x, \xi)= a \cdot x + b \cdot \xi$ with $a$  and $b$
in $H$. From (\ref{2.22}),
$$ ( H_{h/2} F L) (X) = ( H_{h/2} F ) (X) L(X) +  \int _{B^2}  \widetilde F( X +Y )
( \ell _a (y) + \ell _b (\eta) )  d\mu _{B^2, h/2} (Y).  $$
Integrating by parts,
$$ ( H_{h/2} F L) (X) = ( H_{h/2} F ) (X) L(X) +{h\over 2}   \int _{B^2}  \widetilde G( X +Y )
 d\mu _{B^2, h/2} (Y)  $$
where $G(x, \xi) = \Big ( a \cdot \partial _x + b \cdot \partial _{\xi} \Big ) F$.
In other words,
$$ ( H_{h/2} F L) (X) = ( H_{h/2} F ) (X) L(X) +{h\over 2}  \Big ( a \cdot \partial _x
+ b \cdot \partial _{\xi} \Big ) ( H_{h/2} F ) (X). $$
Since $H_{h/2}$ leaves $F$ invariant, this may be written as,
$$  ( H_{h/2} F ) ( H_{h/2} L) +  {h\over 2} \sum _{j\in \Gamma}
\left [ {\partial H_{h/2} F  \over  dx_j} {\partial H_{h/2} L  \over  dx_j} +
{\partial H_{h/2} F  \over  d\xi_j} {\partial H_{h/2} L  \over  d\xi_j} \right ].
$$
Similarly,
$$ H_{h/2} \{ F ,  L \} = \{ H_{h/2} F ,   L \} =   \{ H_{h/2} F ,  H_{h/2} L \}.$$
Consequently, if $\Phi$ is defined in (\ref{b1}) then
$$ H_{h/2} \Phi = ( H_{h/2} F) ( H_{h/2} L) +\hskip 4cm$$
$$\hskip 4cm {h\over 2} \sum _{j\in \Gamma}  \left ( {\partial \over \partial x_j} - i
{\partial \over \partial \xi_j}\right )
( H_{h/2} F)
 \left ( {\partial \over \partial x_j} + i
{\partial \over \partial \xi_j}\right )
( H_{h/2} L).  $$
From Theorem \ref{t2.5}, $ H_{h/2} \Phi$ is the Wick symbol of the composition of the two operators with  Wick symbols being $H_{h/2} F$ and $H_{h/2} L$, that is to say,
$Op_h^{weyl} (F)$ and $Op_h^{weyl} (G)$.
The proposition  is then a consequence of the following Lemma.

\hfill \carre

\begin{lemm}\label{l2.7} Two continuous quadratic forms  on ${\cal D}$ with the same Wick symbol are equal.
\end{lemm}

{\it Proof.} Let $A$ be a continuous quadratic form on ${\cal D}$
which Wick symbol vanishes identically. Let
$X$ and $Y$ be in $H^2$. Set,
$$ \varphi ( \lambda , \mu ) = S_h (A) \left ( {X+Y \over 2} + \lambda {X-Y \over 2} ,
{X+Y \over 2} + \mu {X-Y \over 2} \right ). $$
This function on ${\bf C}^2$ is holomorphic in $\lambda $, anti-holomorphic in
 $\mu$, and identically vanishing if $\lambda = \mu$. It is then identically vanishing
 and the equality $\varphi (1 , -1)=0$ shows that $S_h (A) (X , Y)=0$. The bi-symbol of
 $A$ is identically vanishing.  Let $f$ and $g$ in ${\cal D}_E$ where
 $E \subset H$ is a subspace of finite dimension $n$. Let $C$ and $m$
 be the constants such that we have (\ref{2.8}) for all $f$ and $g$ in ${\cal D}_E$.  Denote by $D(E,m)$ of functions $f$ such that the integral $I(E,m)(f)$ is finite, where $I(E,m)(f)$ is given in (\ref{2.9}). We also have
 $$ f = ( 2\pi h)^{-n} \int _{E^2} < f , \Psi_{Xh} > \Psi_{Xh} dX $$
 and similarly for $g$. Then applying \cite{Y} (Section V.5) one obtains $A(f,g)$ vanishes.

 \hfill\carre

\subsection{Unbounded operators. Sobolev spaces.}\label{s2.F}

We denote by  $W$ the completion of ${\cal  D}$  for the following norm,
$$ \Vert u\Vert_W ^2 = \Vert u \Vert ^ 2 + \sum _{j\in \Gamma }
\Vert ( Q_h(e_j ) + i P_h(e_j )) u\Vert ^2.  $$
Using annihilation operators, one has
$Q_h(e_j ) + i P_h(e_j ) = \sqrt {2h} a_h (e_j)$. Using the number operator
$N = \sum a_h ^{\star } (e_j) a_h (e_j)$, one has  $ \Vert u\Vert_W ^2
= < (I + 2h N) u , u>$ (See also \cite{KR} and \cite{LA2} for other Sobolev spaces in infinite dimension).

\begin{prop}\label{p2.8.}  i) For all $(a, b)$ in $H^2$, let $F_{a , b} (q, p) =
 a \cdot q + b \cdot p$. Then  the operator $Op_h^{weyl} (F_{a b})$ from ${\cal D}$
 into itself,  may be extended to an operator from
$W$ in $L^2(B , \mu _{B , h/2})$ and we have,
\be\label{2.?} \Vert Op_h^{weyl} (F_{a b}) u \Vert  \leq C (|a| + |b|) \ \Vert u \Vert _W.\ee
 ii) Let $F$ in $ S_3(M , \varepsilon )$.
Then the operator $A_h = Op_h^{weyl} (F )$ is bounded from $W $ into $W$.
\end{prop}

{\it  Proof.} i) Point i)  follows from estimates  in Derezi\'nski-G\'erard \cite{DG} , Lemma 2.1 or Lemma 2.3. The operator $Op_h^{weyl} (F_{a b})$
is then denoted by $\Phi_S (a+i b)$.

 ii) For all $u$ in $W$ and for all $j$ in $\Gamma$,
we have from Proposition \ref{p2.6},
$$  ( Q_h(e_j ) + i P_h(e_j )) A_h u =  A_h  ( Q_h(e_j ) + i P_h(e_j ))u + h
Op_h^{weyl } (G_j) u $$
with $G_j (x , \xi) = {\partial F \over \partial x_j} + i {\partial F \over \partial \xi_j}$.
This function belongs to a set $S_2 ( M \varepsilon_j , \varepsilon)$.
From  Theorem 1.4 of \cite{AJN},
$$\Vert A_h \Vert \leq M',\qquad\Vert Op_h^{weyl } (G_j) \Vert \leq M' \varepsilon _j$$
where $M'$ is independent of $j$. The proposition then follows.

\hfill\carre

\section{Reduction to finite dimension.}\label{s3}

With a given bounded operator $A$ in
$L^2 (B, \mu _{B , h/2})$, one always may associate
 a Wick symbol
$\sigma_h^{wick}(A)$. If $A$  verifies the hypotheses of Theorem \ref{t1.2},
we shall associate a Weyl symbol $F$ (which will depend on
$h$). Functions $F$ will satisfy $H_{h/2} F = \sigma_h^{wick}(A)$.

We bring this study to issues related to subspaces $E$ of
finite dimension in $B' \subset H$. One associates two partial heat operators with each
subspace $E \subset B'$. For
any bounded continuous function $F$  on $H^2$ and for all $t > 0$, one set,
\be\label{3.1}  (H_{E , t}  F )(X) = \int _{E^2} F (X + Y_E )
d\mu _{ E^2 , t } (Y_E ).\ee
One can also define a partial heat operator acting, not on
the variables of $E^2$, but on those of its orthogonal. The notation
$E^{\perp}$  now denotes,
\be\label{3.2}  E^{\perp}  = \{ x \in  B, \ \ \  u(x) = 0  \ \ \ \  u \in  E  \}.\ee
This heat operator related to the variables of $(E^{\perp})^2$ can only act on bounded continuous functions $F$ on $H^2$ with
a stochastic extension $\widetilde F$ (bounded measurable function on $B^2$).
One set
\be\label{3.3}   (H_{ E^{\perp} , t } F )(X) = \int _{( E^{\perp})^2}
\widetilde  F (X + Y_{E^{\perp}}  ) d\mu _{ (E^{\perp})^2 , t } (Y_{E^{\perp}}  ).\ee
Indeed, we know from Ramer \cite{RA} (Section 1.B),
that the space $E^{\perp}$ defined in (\ref{3.2}) is also endowed with a gaussian measure. Similarly to $H_t$, we note that,
\be\label{3.4}   \sup_{X\in H^2} | (H_{ E^{\perp} , t } F )(X) |  \leq
\sup_{X\in H^2} |F(X)|.\ee
If $F$ is bounded and continuous on $H^2$ and if its
stochastic extension $\widetilde F$ exits, then we have, from \cite{RA} (Section 1.B,),
\be\label{3.5}  H_{h/2} F = H_{E , h/2 }   H_{E^{\perp}  , h/2 } F.\ee
We then consider an increasing sequence  $(\Lambda _n)$ of finite subspaces in  $\Gamma$ whose union
is $\Gamma$.  We set,
$$E(\Lambda _ n) = {\rm  Vect} (e_j,  \  j \in \Lambda _ n).$$
In Sections \ref{s4} to \ref{s7}, we shall prove the following propositions.

\begin{prop}\label{p3.1}
Let $A$ be a bounded operator in $L^2 (B, \mu _{B , h/2})$ satisfying the
hypotheses of Theorem \ref{t1.2}. Then,

i) the function $\sigma _h^{wick} (A)$ is
in the set $S_{m+4}(M, \varepsilon )$.

 ii) Setting,
 $$ P_{E(\Lambda _n)} (x, \xi ) = \left ( \sum _{j\in \Lambda _n }
e_j (x)e_j,  \sum _{k\in \Lambda _n }  e_k( \xi )e_k \right ),\qquad  (x, \xi ) \in B^2  $$
 and by denoting $ \Vert \cdot \Vert _{\infty } $ the supremum norm  on $H^2$, we have,
\be\label{a10} \Vert \sigma_h^{wick}  (A) -  \sigma_h^{wick}  (A) \circ P_{E(\Lambda _n)}
\Vert _{\infty} \leq  2M \sum _{j\notin \Lambda _n } \varepsilon_j. \ee
\end{prop}

\begin{prop}\label{p3.2}  Let $A$ be a bounded operator in
$L^2 (B, \mu _{B , h/2})$ satisfying the hypotheses in Theorem \ref{t1.2}.
 Then,
for all $n$, there exists a continuous bounded function $F_n$  on $H^2$
such that, if $0 < h < 1$,

i) We have
\be\label{a4}   H_{ E(\Lambda _n), h/2} F_n = \sigma_h^{wick}  (A).\ee

ii) The function $F_n$  is in $S_m ( M_n, \varepsilon)$ with
\be\label{a5}   M_n = M \prod _{j\in \Lambda _n} (1 + K S_{\varepsilon } ^2
h \varepsilon _j ^2 )\ee
where $K$ is a numerical constant and $S_{\varepsilon} $ is defined
in (\ref{1.13}).

iii) If $n < p$ then the function $F_n - F_p$ is in $S_m ( M_{np}, \varepsilon)$
where
\be\label{a6}  M_{np}  = M  \left [ \sum _{j\in \Lambda _p \setminus \Lambda _n}
K (1 + hS_{\varepsilon} ^2 )^2  h \varepsilon _j^2 \right ]
\prod _{j\in \Lambda _p}  (1 + K S_{\varepsilon} ^2 h\varepsilon _j^2  ).\ee
\end{prop}

These propositions will be proved  in Sections \ref{s4} to \ref{s7}. Let us verify that
Theorem \ref{t1.2} follows from these propositions. From  Proposition \ref{p3.2}, the sequence $(F_n)$ converges
to a function $F$ in $S_m ( M' , \varepsilon)$ where $M'$ is defined
in (\ref{1.16}). Let us show that $H_{h/2} F = \sigma_h^{wick} (A)$.
 From Proposition 8.4 in \cite{AJN},
the functions $F_n$ have stochastic extensions $\widetilde F_n$. Then, we may apply the operator $H_{ E(\Lambda _n )^{\perp} , h/2}$ to both sides of equality   (\ref{a4}). We obtain from  (\ref{a4}) and  (\ref{3.5}),
\be\label{a7} H_{h/2} F_n =  H_{ E(\Lambda_n)^{\perp} , h/2}   \sigma_h^{wick} (A).\ee
Let us now take the limit as  $n$ goes to infinity. We have from the point iii) of Proposition \ref{p3.2},
$$| F_n (X)  - F (X)| \leq  M \left [ \sum _{j\notin \Lambda _n} K h
\varepsilon_j^2    \right ]  \prod _{j\in \Gamma }  (1 + K h\varepsilon_j^2 ). $$
 From (\ref{3.4})
we see that, in the sense of the uniform convergence,
\be\label{a8} \lim _{n\rightarrow \infty }  H_{h/2} F_n =
H_{h/2} F.\ee
We shall also check that,
\be\label{a9}  \lim _{n\rightarrow \infty }    H_{ E(\Lambda _n)^{\perp} , h/2}
  \sigma_h^{wick} (A)  =  \sigma_h^{wick} (A).\ee
Indeed, setting,
$\Psi  = \sigma_h^{wick} (A)$,
we have
$$  \Vert \Psi  - H_{ E(\Lambda _n)^{\perp} , h/2} \Psi \Vert _{\infty }
\leq \Vert \Psi  -  \Psi   \circ P_{E(\Lambda _n ) } \Vert _{\infty } +
\Vert  H_{ E(\Lambda _n)^{\perp} , h/2} ( \Psi  -  \Psi   \circ P_{E(\Lambda _n ) } )
\Vert _{\infty }.$$
We have used the fact that $H_{ E(\Lambda_n)^{\perp} , h/2} (\Psi \circ  P_{E(\Lambda _n)})
= \Psi \circ  P_{E(\Lambda _n)}$. The
limit in (\ref{a9}) follows from (\ref{3.4})(\ref{a10}) and of point  ii) in Proposition \ref{p3.1}. Using (\ref{a7})(\ref{a8})(\ref{a9}) we obtain
$H_{h/2}F = \sigma_h^{wick} (A)$.  Since the function $F$ is in
$S_m(M', K \varepsilon )$ then a Weyl quadratic form is associated with,  by Theorem \ref{t2.2}, and
a bounded operator $ Op_h^{weyl} (F )$ associated with, by Theorem 1.4 of \cite{AJN}.
From (\ref{2.23}), the Wick symbol of this operator is $H_{h/2}F$.
Consequently the operators $Op_h^{weyl} (F )$ and $A$ have the same Wick symbol. From Lemma \ref{l2.7}, these two  operators are equal.

Once Propositions \ref{p3.1} and \ref{p3.2} proved, we have indeed found a function $F$ in $S_m(M', K \varepsilon )$ whose corresponding Weyl operator
equals to $ A$. Theorem \ref{t1.2} is then a consequence of Propositions \ref{p3.1} and \ref{p3.2}.

\hfill\carre

\section{Proof of  Proposition \ref{p3.1}.}\label{s4}

Let $A$ be a bounded operator $A$  in
$L^2 (B, \mu _{B , h/2})$ satisfying the hypotheses of
Theorem \ref{t1.2}. From  Theorem \ref{t2.5}, we have,
\be\label{4.1} \sigma_h^{wick} (  [Q_h(ej ), A] ) =
 ih {\partial \over \partial \xi _j }  \sigma_h^{wick}  (A),\qquad
 \sigma_h^{wick} (  [P_h(ej ), A] )    =
 - ih {\partial \over \partial \xi _j }  \sigma_h^{wick} (A).\ee
For all bounded operator $B$, one has,
$| \sigma_h^{wick} (B)(X)| \leq  \Vert B \Vert $. Consequently,
if $A$ verifies the hypotheses of  Theorem \ref{t1.2} one deduces estimates, for each multi-index
$(\alpha, \beta)$ in  ${\cal M} _{m+4}$,
$$ | \partial_x^{\alpha}\partial_{\xi}^{\beta}  \sigma_h^{wick}  (A)(x, \xi)|
\leq  M \prod _{j\in \Gamma } \varepsilon_j ^{\alpha _j + \beta _j},  $$
which prove point i) of Proposition \ref{p3.1}. We deduce,
$$ | \sigma_h^{wick}(A)(x, \xi )
-  \sigma_h^{wick} (A)(P_{E(\Lambda _n)} (x, \xi)) | \leq 2 M \sum _{j\notin \Lambda _n  }
\varepsilon_j ^{\alpha _j + \beta _j}, $$
which proves Proposition \ref{p3.1}. We shall also need analogous estimates on the bi-symbol. One deduces from (\ref{4.1}) these estimates by setting,  for all $j\in \Gamma$,
 $$ {\partial \over \partial X_j}   = {1\over 2}  \left (  {\partial \over \partial x_j}
 - i  {\partial \over \partial \xi _j}        \right ),\qquad
{\partial \over \partial \overline Y_j}   = {1\over 2}  \left (  {\partial \over \partial y_j}
 + i  {\partial \over \partial \eta _j}        \right ).
            $$
With these notations, one has,
\be\label{4.2}
S_h([Q_h(e_j ), A])(X, Y ) = - h \left (  {\partial \over \partial X_j}  - {\partial \over \partial \overline Y_j}
   \right )   (S_hA)(X, Y ),\ee
\be\label{4.3} S_h([P_h(e_j ), A])(X, Y ) = -i h \left (  {\partial \over \partial X_j}  + {\partial \over \partial \overline Y_j}
   \right )   (S_hA)(X, Y ).\ee
 Consequently, for all multi-indices $(\alpha, \beta )$,
\be\label{4.4}   S_h ( ({\rm ad} P_h)^{\alpha}  ({\rm ad} Q_h)^{\beta}     A )(X, Y )
 = c_{\alpha \beta}
  h^{|\alpha + \beta |} (\partial _x + \partial _y)^{\alpha}  (\partial _{\xi} + \partial _{\eta} )^{\beta }
  S_h(A)(X, Y )\ee
where $|c_{\alpha \beta} | =1$.  With (\ref{2.4}), we deduce that
\be\label{4.5} | (\partial _x + \partial _y)^{\alpha}  (\partial _{\xi} + \partial _{\eta} )^{\beta }
  S_h(A)(X, Y ) | \leq  h^{-|\alpha + \beta |}
  e^{-{1\over 4h} |X-Y|^2} \Vert ({\rm ad} P_h)^{\alpha}  ({\rm ad} Q_h)^{\beta}     A \Vert.\ee

\section{Finite dimensional analysis.}\label{s5}

We consider here the case where $H$ is a real Hilbert space  with finite dimension $n$. Let $A$ be an operator satisfying
 hypothesis of Theorem \ref{t1.2}. Let $\Phi = S_hA$  its bi-symbol,
defined in (\ref{2.11}). We have seen that $\Phi (X , Y)$ is holomorphic in $X$, anti-holomorphic in $Y$.
From (\ref{4.5}), the following norm is finite,
\be\label{5.1} N_h^{(2)} (\Phi) = \sum _{(\alpha , \beta) \in {\cal M} _2}
\Vert e^{-{1 \over 4h} |X - Y |^2 } (\partial _x + \partial _y)^{\alpha }
 (\partial _{\xi } + \partial _{\eta })^{\beta } \Phi  \Vert_{\infty },
\ee
 where $\Vert \cdot \Vert _{\infty }$ is the supremum norm. Note again that a choice of particular basis has been made.

One introduces in  distributions sense an
 integral transform giving  the Weyl symbol $F$ of $A$ starting from the
 bi-symbol $\Phi$, and give estimates on $F$. This integral  is not converging but has to be understood as an oscillatory integrals
 (see H\"ormander \cite{HO}).  This leads to a proof of Beals's theorem in finite dimension
 (see Unterberger \cite{U}). Setting,
\be\label{5.2} K_h^{Beals} (X , Y, Z) =   e^{- {1\over h}(Z-Y)\cdot (\overline Z- \overline X) - {1\over 2h}|X-Y|^2 }.\ee

\begin{theo}\label{t5.1} Let  $H$  be a real Hilbert space of
 finite dimension $n$.  Set $(X , Y) \rightarrow \Phi (X , Y)$ a function
 on $H^2 \times H^2$ which is holomorphic in $X$ and anti-holomorphic in $Y$, such that the norm
 $ N_h^{(2)} (\Phi)$ defined in (\ref{5.1}) is finite (for some orthonormal basis).
 Then,

 i) The following integral transform defines, a priori in the sense of  distributions,
 a function $B _h \Phi $ which is bounded and continuous on $H^2$,
\be\label{5.3} (B_h \Phi ) (Z) = 2^n (2 \pi h)^{-2n} \int _{H^4}
 \Phi (X , Y)  K_h^{Beals} (X , Y, Z) dX dY.\ee
Moreover,  this function satisfies,
\be\label{5.4} \Vert B_h \Phi \Vert _{\infty } \leq K^n N_h^{(2)} (\Phi) \ee

 ii) Moreover, one has,
 \be\label{5.5} ( H_{h/2} B_h \Phi ) (Z) = \Phi (Z , Z).\ee
\end{theo}

 {\it Proof of i).} We follow the method of Unterberger \cite{U}.  The change of variables
 $$ X = Z + S + {T\over 2},\qquad Y = Z + S -  {T\over 2} $$
allows  to rewrite (\ref{5.3}) as,
\be\label{5.6} (B_h \Phi )(Z) = 2^n (2\pi h)^{-2n} \int _{H^4} \Psi (S, T , Z)
K_h(S, T )dSdT\ee
 with
\be\label{5.7} \Psi (S, T , Z ) = \Phi \left  ( Z + S +{T\over 2}  , Z + S - {T\over 2} \right  )
\ee
\be\label{5.8}  K_h(S, T ) = e^{-{1\over h} |S|^2 - {i \over h}\sigma (S , T) - {1\over 4h} |T|^2}.\ee
Set $S_j = (s_j, \sigma _j )$, $T_j = (t_j, \tau _j)$. Let $L_j$ and $M_j$ be the
operators defined, for each function $G(S, T )$, by
$$ L_j G = \left ( 1 + {\tau _j ^2 \over h } \right ) ^{-1} e^{-{1\over h} s_j^2 }
\left (  1 - h {\partial ^2 \over \partial s_j^2 }      \right ) e^{{1\over h} s_j^2 } G  $$
$$ M_j G =  \left ( 1 + {t _j ^2 \over h } \right ) ^{-1}  e^{-{1\over h} \sigma _j^2 }
\left (  1 - h {\partial ^2 \over \partial \sigma _j^2 }      \right )  e^{{1\over h} \sigma _j^2 } G. $$
One verifies that,
$$ L_jK_h = K_h,\qquad  M_jK_h = K_h \hskip 2cm  j\leq n   $$
where the function $K_h$ defined in (\ref{5.8}).
 Consequently,
 $$  (B_h \Phi )(Z) = 2^n (2\pi h)^{-2n} \int_{H^4}  K_h(S, T )
 \left [ \prod _{j\leq n}   ^tL_j \ ^ t M_j \right ] \Psi  	(S, T , Z )
 dS dT. $$
We see that,
$$ ^tL_j =  \left ( 1 + {\tau _j ^2 \over h } \right ) ^{-1} \Big [
  a_0(s_j / \sqrt h) + h^{1/2}  a_1(s_j / \sqrt h)\partial _{s_j} +
  h a_2(s_j / \sqrt  h) \partial _{s_j} ^2 \Big ]  $$
with
$$  a_0(s) = 3 - 4s^2,\qquad
a_1(s) = 4s,\qquad
a_2(s) = -1.$$
Similarly,
 $$ ^ tM_j =  \left ( 1 + {t _j ^2 \over h } \right ) ^{-1}
 \Big [
  a_0(\sigma _j / \sqrt h) + h^{1/2}  a_1(\sigma _j / \sqrt h)\partial _{\sigma_j} +
  h a_2(\sigma _j / \sqrt  h) \partial _{\sigma _j} ^2 \Big ] $$
Consequently,
$$  | (B_h  \Phi )(Z)| \leq  \sum _{ (\alpha , \beta) \in {\cal M}_2 }
h^{|\alpha + \beta |/2} F_{\alpha \beta } (Z) $$
$$  F_{\alpha \beta } (Z) =  2^n (2\pi h)^{-2n} \int _{H^4}
 e ^{ -{1\over h} |S|^2} \prod _{j\leq n} \left ( 1 + {t _j ^2 \over h } \right ) ^{-1}
  \left ( 1 + {\tau _j ^2 \over h } \right ) ^{-1} \left | a^{\alpha } (s / \sqrt h)
   a^{\beta } (\sigma  / \sqrt h)        \right |   $$
$$ | e ^{ -{1\over 4h} |T|^2}  \partial _s ^{\alpha }   \partial _{\sigma} ^{\beta  }
   \Psi 	(S, T , Z )|   dS dT $$
 where we have set
 $$  a ^{\alpha } (s ) = \prod _{j\leq n} a_{\alpha _j}  (s_j ).  $$
 There exists $K > 0$ such that,
$$ \pi ^{-1/2} \int _{\R}e^{-s^2} |a_j(s)| ds \leq  K,\qquad  0 \leq  j \leq  2 $$
 and also
 $$ (2 \pi ) ^{-1/2}  \int _{\R}  (1 + x^ 2)^{-1} dx \leq  K. $$
 Consequently,
 $$  | (B_h  \Phi )(Z)|  \leq K^n \sum _{ (\alpha , \beta) \in {\cal M}_2 }
h^{|\alpha + \beta |/2}  \sup _{(S , T) \in H^4}
\left |  e^{-{1\over 4h} |T|^2}  \partial_s ^{\alpha }   \partial_{\sigma }  ^{\beta  }
\Psi (S , T, Z)       \right |.  $$
From the defintion of $\Psi $ in (\ref{5.7}),
$$  | (B_h  \Phi )(Z)|  \leq K^n  \sum _{ (\alpha , \beta) \in {\cal M}_2 }
h^{|\alpha + \beta |/2}  \sup _{(X , Y) \in H^4}
\left |  e^{-{1\over 4h} |X-Y|^2 } ( \partial_x + \partial_y) ^{\alpha }
( \partial_{\xi } + \partial_{\eta} ) ^{\beta  }
\Phi (X  Y)       \right |.$$
We then deduce (\ref{5.5}) with another constant $K$.

 {\it Proof of ii).} If a function $\Psi$ on $H^2$ is written as
 $$ \Psi (Z) = e^{ -{1\over h} |Z|^2 + {1\over h} (A \cdot Z + B \cdot \overline Z )}$$
 where $A$ and $B$ are in $H^2$, and $A \cdot Z $ denotes the bi-{\bf C}-linear scalar product, then the action of  the heat operator on $\Psi$
 verifies,
 $$ ( H_{h/2} \Psi ) (Z) = \Big ( e^{{h\over 4} \Delta }  \Psi  \Big )
 (Z)  = 2^{-n} e^{{1\over 2h} |Z|^2 + {1\over 2 h} (A \cdot Z + B \cdot \overline Z )
   {1\over 2 h} A \cdot B }.$$
Thus,
$$ (  H_{h/2} K_h^{Beals} (X , Y, \cdot ) ) (Z) = 2^{-n}
{\cal B}_h (Z , Z , X , Y) e^{-{1\over 2h} (|X|^2 + |Y|^2)} $$
where ${\cal B}_h$ is our type of reproducing kernel introduced in (\ref{2.19}). Consequently, as in (\ref{2.20})
$$  ( H_{h/2} B_h \Phi ) (Z) = \int _{H^4} \Phi (X , Y) {\cal B}_h (Z , Z , X , Y)
d\mu _{H^4 , h } (X , Y) = \Phi (Z , Z).$$

\section{Proof of Proposition 3.2: first step. }\label{s6}

For all operators $A$ satisfying the
hypotheses of Theorem \ref{t1.2} and for some subsets $E $ of
finite dimension in $B' \subset H$, we shall find  a bounded  continuous function
$\tau_{E,h} (A)$   on $H^2$ such that
\be\label{6.1} H_{E, h/2}   \tau_{E,h} (A) = \sigma_h^{wick}  (A).\ee
It is point i) of Proposition \ref{p3.2}. Moreover, we shall give estimations on  this function. For all finite
subsets $I$ in $\Gamma $, let $E(I)$ be the subspace of $B' \subset H$
spanned by the $e_j $, $j \in I$. Recall that the elements $e_j$
$(j \in \Gamma)$ of our Hilbertian basis are in $B'$. Let
${\cal  M}_2(I)$ be the set of all multi-indices $(\alpha, \beta)$ such that
$\alpha _j = \beta _j = 0$  if $j  \notin I$, and $\alpha _j \leq  2$
and $\beta _j \leq  2$  if $j \in I$.

\begin{prop}\label{p6.1} Let $A$ be an operator  satisfying the hypotheses in Theorem \ref{t1.2}. Set $I$
a finite subspace of $\Gamma$.  Then, there exists a bounded continuous function
$ \tau_{E(I),h} (A)$  on $H^2$ satisfying (\ref{6.1}). Moreover,
\be\label{6.2} \Vert     \tau_{E(I),h} (A) \Vert _{\infty }  \leq
K^{|I|}  \sum _{(\alpha , \beta)\in {\cal M}_2(I) }  h^{-|\alpha +\beta|/2}
\Vert ({\rm ad} P_h)^{\alpha}   ({\rm ad} Q_h)^{\beta} A \Vert
\ee
 where  $ K$  is a numerical constant.
 \end{prop}

{\it Proof.} We denote  $E = E(I)$, $E^{\perp}$ the orthogonal complement of $E$
in $H$, and $Z = (Z_E, Z_{E^{\perp}} )$ the variable in $H^2$.
For all $Z_{E^{\perp}} $ in $(E^{\perp})^2$, we shall apply Proposition \ref{t5.1} replacing $H$ by $E$, with the following function
$\Phi $ defined on $E^2$,
$$ \Phi_ { Z_{E^{\perp}} } (X_E , Y_E) = (S_hA) ( X_E , Z_{E^{\perp}} , Y_E,
Z_{E^{\perp}} ). $$
Using again notation (\ref{5.3}), which a priori only makes sense as an oscillatory integral on $E^2$, one set for all  $Z = (Z_E, Z_{E^{\perp}} )$ in $H^2$,
$$  \tau_{E(I),h} (A) (Z) =  2^{ {\rm dim} (E)}  (2 \pi h)^{-2{\rm dim} (E) }
\int _{E^4}
(S_hA) ( X_E , Z_{E^{\perp}} , Y_E,
Z_{E^{\perp}} )   K_h^{Beals} (X_E , Y_E, Z_E) dX_E dY_E$$
where $ K_h^{Beals} $ is defined in (\ref{5.2}). One may apply  Theorem \ref{t5.1}, choosing as an orthonormal basis  of $E = E(I)$, the one constituted with the $e_j$ $j\in I$.
With this choice, we have from (\ref{4.5}),
$$ N_h^{(2)} (\Phi_ { Z_{E^{\perp}} } ) \leq  \sum _{(\alpha , \beta)\in {\cal M}_2(I) }  h^{-|\alpha +\beta|/2}
\Vert ({\rm ad} P_h)^{\alpha}   ({\rm ad} Q_h)^{\beta} A \Vert $$
and the term in the right hand side is finite under  hypothesis of Theorem \ref{t1.2}.
From  Theorem \ref{t5.1}, the function $  \tau_{E(I),h} (A)$ is well-defined,
continuous and bounded on $H^2$ and satisfies (\ref{6.1}) and (\ref{6.2}).

\hfill\carre

\section{Proof of Proposition \ref{p3.2}: second step.}\label{s7}

For all finite subsets  $I$ of $\Gamma$, let us set
\be\label{7.1}  T_{I, h}  =  \prod _{j\in I} (I - H_{D_j , h/2} )\ee
where $D_j$ is spanned by the vector $e_j$  of our Hilbertian basis of $H$, and $H_{D_j , h/2}$  is the operator defined in (\ref{3.1}),
with $E$ replaced by $D_j$, thus with an integral on $D_j ^2 $.
When $I =\emptyset$, we set  $T _{I , h} = I d$.
We denote by $E(I)$ the subspace
of $B'$  spanned by the $e_j$, $j\in I$. Recall that the elements $e_j$
 $(j \in \Gamma)$   of our Hilbertian basis of $H$ are in $B'$.
If $I = \emptyset$ then
set $E(I) = \{ 0 \}$. For any operator $A$ satisfying the
hypotheses in Theorem \ref{t1.2} and for all subspaces $E \subset  B' \subset  H$
of finite dimension, set $\tau _{E(I) , h} (A)$ the function on $H^2$ defined in
the Proposition \ref{p6.1}. In particular, we may have $E = E(I)$ with $I$
being a finite subset of $\Gamma$.  We choose an increasing sequence
$(\Lambda_n)$ of finite subsets of $\Gamma$ with its union equals to $\Gamma$.
 For all $ n$, one defines a function $F_n $ on $H^2$ by,
\be\label{7.2}   F_n = \sum _{I\subset \Lambda _n} T _{I , h}
 \tau_{E(I) , h} (A).\ee
The above sum is running over all the subsets $I$ of $\Lambda _n$  including the empty set. We shall show  that this sequence
of functions has indeed the properties announced
Proposition \ref{p3.2}.

{\it  Point i) }
One has, for all subsets $I \subset  \Lambda _n$,
$$ H_{E(\Lambda _n), h/2}  =  H_{E(I), h/2}    H_{E(\Lambda _n \setminus I), h/2}   $$
 and these operators  commutes with each other and with $T _{I , h}$.
 Consequently,
 $$ H_{E(\Lambda _n), h/2} F_n  =  \sum _{I\subset \Lambda _n}  T _{I, h}
 H_{E(\Lambda _n \setminus I), h/2}   H_{E(I), h/2}  \tau_{E(I) , h} (A). $$
 From equality (\ref{6.1})
applied to set $E(I)$, one has,
$$ H_{E(\Lambda _n), h/2} F_n  =  \sum _{I\subset \Lambda _n}    T _{I , h}
 H_{E(\Lambda _n \setminus I), h/2} \sigma_h^{wick} (A). $$
The following equality is
a variant of the binomial formula,
$$\sum _{I\subset \Lambda _n}  T _{I , h}  H_{E(\Lambda _n \setminus I), h/2} = I d.
$$
So, we have proved equality (\ref{a4}), point i) of the Proposition \ref{p3.2}.

Points ii) and iii)
will both be a direct consequence of the following inequality. If $A$
satisfies hypothesis in Theorem \ref{t1.2}, for all $(\alpha , \beta)$
 in $M_m$, for any finite subset $I$ in $\Gamma$ and for all
  $h$ in $(0, 1)$,
\be\label{7.3}  \Vert \partial_z ^{\alpha }  \partial_{\zeta}  ^{\beta } T_{I, h}
  \tau_{E(I) , h} (A) \Vert _{\infty} \leq  M (K S_{\varepsilon} ^ 2) ^{|I|}
  \prod _{j\in I}  h \varepsilon_j^2 \prod _{j\in \Gamma } \varepsilon _j ^{\alpha _j + \beta _j}
\ee
 where $K$ is a numerical constant  and $S_{\varepsilon}$  is
defined in (\ref{1.13}).

It remains to prove (\ref{7.3}). If $H_{D_j,h/2}$ is defined in (\ref{3.1}),
with $E$  replaced by
$D_j = {\rm Vect}\,(ej )$, we may write,
$$  I - H_{ D_j, h/2}  = {h\over 4}  V_j
(\partial _{z_j} ^2  + \partial _{\zeta _j} ^2 )  $$
where the operators $V_j$ are bounded in the space $ C_b$ of
continuous bounded functions on $H^2$, and are commuting with partial derivatives operators. Moreover,
$$ \Vert  V_j \Vert _{{\cal L} (C_b)}  \leq  1. $$
 Therefore, one may rewrite the operator $T_{I, h}$  defined in
 (\ref{7.1}) under the following form,
$$ T_{I, h} = \prod _{j\in I}  (h/4)  V_j (\partial _{z_j} ^2  + \partial _{\zeta _j} ^2 ).$$
 Let $ {\cal N} (I) $ be the set of
multi-indices $(\alpha , \beta )$ such that $\alpha _j = \beta _j = 0$ if
$ j \notin  I$,  and if $j \in  I$,
either we have $\alpha _ j = 2$  and $\beta _j = 0$, or
$\alpha _j = 0$ and $\beta _j = 2$. Consequently,
$$\Vert \partial_z ^{\alpha }  \partial_{\zeta}  ^{\beta } T_{I, h}
  \tau_{E(I) , h} (A) \Vert _{\infty} \leq (h/4)^{|I|} \sum _{ (\gamma , \delta) \in {\cal N} (I) }
 \Vert \partial_z ^{\alpha + \gamma  }  \partial_{\zeta}  ^{\beta + \delta  }
  \tau_{E(I) , h} (A) \Vert _{\infty}.$$
On verifies that,
$$    \left [  {\partial \over \partial x_j} +  {\partial \over \partial y_j} + {\partial \over \partial z_j}
     \right ] K_h^{Beals} (X, Y, Z ) = 0,\qquad
     \left [  {\partial \over \partial {\xi}_j} +  {\partial \over \partial {\eta}_j} + {\partial \over \partial \zeta _j}
     \right ] K_h^{Beals} (X, Y, Z ) = 0.
 $$
Consequently,
$$  \partial_z ^{\alpha }  \partial_{\zeta}  ^{\beta } \tau _{E(I),h} (A) =
\tau _{E(I),h} A_{\alpha \beta}  $$
 where $A_{\alpha \beta}$  is such that,
 $$ (S_h A_{\alpha \beta} ) (X, Y ) =  ( \partial _x +\partial _y)^{\alpha }
  ( \partial _{\xi} +\partial _{\eta} )^{\beta }  (S_hA)(X, Y ).$$
 From (\ref{4.4}),
 $$A_{\alpha \beta} = c_{\alpha \beta} h^{-|\alpha +\beta |} ({\rm ad} P_h)^{\alpha}
 ({\rm ad} Q_h)^{\beta} A$$
where $|c_{\alpha \beta} |  = 1$. Then,
$$ \Vert   \partial_z ^{\alpha }  \partial_{\zeta}  ^{\beta }  T _{I , h}
   \tau _{E(I),h} (A) \Vert _{\infty }  \leq  (h/4)^{|I|} \sum _{(\gamma , \delta ) \in {\cal N} (I)}
h^{ -|\alpha + \beta + \gamma + \delta | }
\Vert \tau _{E(I),h} ( ({\rm ad} P_h)^{\alpha + \gamma }
 ({\rm ad} Q_h)^{\beta + \delta } A \Vert _{\infty }.
    $$
From the Proposition \ref{p6.1},
$$ \Vert   \partial_z ^{\alpha }  \partial_{\zeta}  ^{\beta }  T _{I , h}
   \tau _{E(I),h} (A) \Vert _{\infty }  \leq  (K h/4)^{|I|}
    \sum _{(\gamma , \delta ) \in {\cal N} (I)}
     \sum _{(\lambda  , \mu ) \in {\cal M}_2 (I)}
     h^{ -|\alpha + \beta + \gamma + \delta | - | \lambda + \mu |/2 }$$
 $$ \hskip 3cm \Vert ( ({\rm ad} P_h)^{\alpha + \gamma + \lambda  }
 ({\rm ad} Q_h)^{\beta + \delta + \mu  } A \Vert. $$
If $(\alpha , \beta) \in {\cal  M} _m$, $(\gamma  , \delta ) \in {\cal  N}  (I)$
  and $(\lambda , \mu) \in {\cal  M} _2(I)$, then the sum
  $ (\alpha  + \gamma  + \lambda , \beta  + \delta  + \mu)$ belongs to
  $ {\cal M} _{m+4}$.
From assumptions of Theorem \ref{t1.2},
$$  \Vert   \partial_z ^{\alpha }  \partial_{\zeta}  ^{\beta }  T _{I , h}
   \tau _{E(I),h} (A) \Vert _{\infty }  \leq
 M(K h/4)^{|I|}  \sum _{(\gamma , \delta ) \in {\cal N} (I)}
     \sum _{(\lambda  , \mu ) \in {\cal M}_2 (I)}
     h^{ |\lambda + \mu | /2} \prod _{j\in \Gamma} \varepsilon _j^{ \alpha _j
     + \beta _j + \gamma _j + \delta _j + \lambda _j + \mu_j}.  $$
The number of multi-indices in ${\cal N} (I)$ is $2^{|I|}$,
and the number of multi-indices in ${\cal M}_2(I)$ is
$9^{|I|}$. For all multi-indices $ (\gamma , \delta) \in {\cal  N}  (I)$, we have
$ \gamma _j + \delta _j = 2$  if $j\in I$.  If $0 < h < 1$, for all multi-indices
$ (\lambda , \mu ) \in {\cal  M}_2(I)$,  we have $ (
\sqrt h  \varepsilon _j )^{\lambda _j+ \mu _j}  \leq  S_{\varepsilon}^ 2$,
where $ S_{\varepsilon}$  is defined in (\ref{1.13}). Consequently, we have indeed proved (\ref{7.3}) with another universal constant $K$.
From (\ref{7.2}), we deduce the points ii) and iii) of the
Proposition \ref{p3.2}, which complete the proof of Theorem \ref{t1.2}.

\hfill\carre

\section{Composition of operators.}\label{s8}

\begin{theo}\label{t8.1}  Let $F$ in
$S_{m+6}(M , \varepsilon )$ and $G$ in $S_{m+6}(M' , \varepsilon )$ ($m \geq  0$).
 Then there exists a function $ H_h$ in $S_{m}(M'' ,(m+4) \varepsilon ) $ such that,
\be\label{8.1}  Op_h^{weyl} (F ) \circ  Op_h^{weyl}(G) = Op_h^{weyl}(H_h).\ee
We have set,
$$ M''   = M M'  \prod _{j\in \Gamma} (1
+ K(m + 4)^2  S_{\varepsilon}^2 h \varepsilon_j ^2)^3  \leqno (8.2) $$
where $K$ is a universal constant and $S_{\varepsilon}$
is defined in (\ref{1.13}).
\end{theo}

{\it Proof.}  For any multi-index $(\alpha, \beta)$
in ${\cal M}_{m+4}$ we have,
$$  ({\rm ad}P_h )^{\alpha} ({\rm ad}Q_h )^{\beta}  \Big ( Op_h^{weyl} (F ) \circ Op_h^{weyl} (G) \Big )
 = $$
$$ \sum_{ \alpha ' + \alpha '' = \alpha \atop \beta ' + \beta '' = \beta }
\Big ( ({\rm ad}P_h )^{\alpha '} ({\rm ad}Q_h )^{\beta' }  Op_h^{weyl} (F )      \Big )   \circ
\Big ( ({\rm ad}P_h )^{\alpha ''} ({\rm ad}Q_h )^{\beta'' }  Op_h^{weyl} (G )      \Big ).
 $$
From (\ref{1.15}) (with $m$ replaced by $m +6$) and similarly
for $G$, we have, for each multi-index  $(\alpha, \beta)$ in
${\cal M}_{m+4}$,
$$ ({\rm ad}P_h )^{\alpha} ({\rm ad}Q_h )^{\beta}
\Big ( Op_h^{weyl} (F ) \circ Op_h^{weyl} (G) \Big ) \Vert  \leq
 M M' N(\alpha , \beta  ) \prod _{j\in \Gamma}
 (1+ 81 \pi h S_{\varepsilon} \varepsilon _j^2)^2
  \prod _{j\in \Gamma} (h\varepsilon _j)^{\alpha _j + \beta _j}
$$
where $N(\alpha, \beta )$ is the number of decompositions of $(\alpha, \beta)$
as a sum  of two multi-indices $(\alpha ', \beta ')$ and $(\alpha '', \beta '')$.
If $(\alpha, \beta)$ is in ${\cal M}_{m+4}$ then this number equals is smaller than $ (m + 4)^{|\alpha + \beta|}$.
 Consequently, $Op_h^{weyl} (F ) \circ Op_h^{weyl} (G)$ satisfies
a condition similar to (\ref{1.15b}) with $\varepsilon _j$
remplaced by $(m + 4)\varepsilon _j$. So our Theorem \ref{t8.1} is
is  a consequence of Theorem \ref{t1.2}.

\hfill\carre

\medskip

laurent.amour@univ-reims.fr\newline
LMR EA 4535 and FR CNRS 3399, Universit\'e de Reims Champagne-Ardenne,
 Moulin de la Housse, BP 1039,
 51687 REIMS Cedex 2, France.

rlascar@math.univ-paris-diderot.fr\newline
Institut Mathématique de Jussieu UMR CNRS 7586,  Analyse Algébrique, 4 Place Jussieu, 75005 Paris, France.

jean.nourrigat@univ-reims.fr\newline
LMR EA 4535 and FR CNRS 3399, Universit\'e de Reims Champagne-Ardenne,
 Moulin de la Housse, BP 1039,
 51687 REIMS Cedex 2, France.

 \end{document}